\documentclass[lean]{AFM}
\usepackage{array}

\renewcommand{\.}{\hskip .75pt}

\newcommand{\fin}[1]{[\.#1\.]}

\newcommand{\sekt}[1]{Section~\ref{#1}}

\DeclareMathOperator{\aand}{\;\wedge\;}
\DeclareMathOperator{\st}{,\;}

\let\r=\rightarrow
\let\*=\cdot

\makeatletter
\newcommand{\thickhline}{%
	\noalign {\ifnum 0=`}\fi \hrule height 2pt
	\futurelet \reserved@a \@xhline
}
\newcolumntype{"}{@{\hskip\tabcolsep\vrule width 2pt\hskip\tabcolsep}}
\makeatother

\makeatletter
\newcommand{\customlabel}[2]{%
   \protected@write \@auxout {}{\string \newlabel {#1}{{#2}{\thepage}{#2}{#1}{}} }%
   \hypertarget{#1}{#2\!\!}
}
\makeatother

\makeatletter
\newcommand{\customlabelrestated}[2]{%
   \protected@write \@auxout {}{\string \newlabel {#1}{{#2}{\thepage}{#2}{#1}{}} }%
   \hypertarget{#1}{#2, restated\!\!}
}
\makeatother

\title[Duality theory in linear optimization and its extensions]{Duality theory in linear optimization and its extensions: formally verified}

\author[M. Dvorak and V. Kolmogorov]{Martin Dvorak and Vladimir Kolmogorov} 

\authorinfo[M. Dvorak]{Institute of Science and Technology Austria}{martin.dvorak@ista.ac.at}

\authorinfo[V. Kolmogorov]{Institute of Science and Technology Austria}{vnk@ista.ac.at}

\VOLUME{2}
\YEAR{2026}
\NUMBER{1}
\firstpage{1}
\DOI{https://doi.org/10.46298/afm.14253}
\receiveddate{September 13, 2024}
\finaldate{July 17, 2025}
\accepteddate{September 23, 2025}

\msc{68V20, 15A39, 90C05}

\begin{abstract}
Farkas established that a system of linear inequalities has a solution if and only if we cannot obtain
a contradiction by taking a linear combination of the inequalities.
We state and formally prove several Farkas-like theorems over linearly ordered fields in Lean~4.
Furthermore, we extend duality theory to the case when some coefficients are allowed to take
``infinite values''.
\end{abstract}

\keywords{Farkas lemma, linear programming, extended reals, calculus of inductive constructions}

\begin{document}

\section{Introduction}
\label{introduction}

In optimization and related fields, we often ask whether a given
system of linear inequalities has a solution.
Duality theorems answer these naturally arising questions as follows:
A system of linear inequalities has a solution if and only if
we cannot obtain a contradiction by taking a linear combination of
the inequalities. When we formulate duality theorems as
``either there exists a solution, or there exists a vector of coefficients
that tells us how to derive $0 < 0$ from given inequalities'',
they are usually called {\em theorems of alternatives}.
Two well-known theorems of alternatives are as follows
(Farkas~\cite{Farkas1894,Farkas1898}; Minkowski~\cite{Minkowski}).

\begin{theorem}[\customlabel{equalityFarkas}{\texttt{equalityFarkas}}]
Let $I, J$ be finite types. 
Let $F$ be a linearly ordered field.
Let $A$ be a matrix of type $(I \!\times\! J) \r F$.
Let $b$ be a vector of type $I \r F$.
Exactly one of the following exists:
\begin{itemize}
\item nonnegative vector $x : J \r F$ such that $A \* x = b$
\item vector $y : I \r F$ such that $A^T\! \* y \ge 0$ and $b \* y < 0$
\end{itemize}
\end{theorem}

\begin{theorem}[\customlabel{inequalityFarkas}{\texttt{inequalityFarkas}}]
Let $I, J$ be finite types. 
Let $F$ be a linearly ordered field.
Let $A$ be a matrix of type $(I \!\times\! J) \r F$.
Let $b$ be a vector of type $I \r F$.
Exactly one of the following exists:
\begin{itemize}
\item nonnegative vector $x : J \r F$ such that $A \* x \le b$
\item nonnegative vector $y : I \r F$ such that $A^T\! \* y \ge 0$ and $b \* y < 0$
\end{itemize}
\end{theorem}

For optimization problems, we obtain a correspondence between
the optimal value of a linear program and the optimal value of its dual
(originally discussed in the context of zero-sum games by Dantzig and von Neumann~\cite{Discussion};
later in Gale, Kuhn, Tucker~\cite{GaleKuhnTucker}).
The strong duality theorem, which is a cornerstone of linear programming,
can be stated as follows.

\begin{theorem}[\customlabel{StandardLP.strongDuality}{\texttt{StandardLP.strongDuality}}]
Let $I$ and $J$ be finite types.
Let $F$ be a linearly ordered field.
Let $A$ be a matrix of type $(I \times J) \r F$.
Let $b$ be a vector of type $I \r F$.
Let $c$ be a vector of type $J \r F$.
Then
\begin{equation}
\min \,\{\, c \* x ~|~ x \ge 0 \aand A \* x \le b \,\}
=- \min \,\{\, b \* y ~|~ y \ge 0 \aand (-A^T) \* y \le c \,\}
\label{eq:StandardLP.strongDuality:minmin}
\end{equation}
holds if at least one of the systems has a solution (very roughly paraphrased).
\end{theorem}
Using identity $\min\{c\*x~|~\ldots\}=-\max\{-c\*x~|~\ldots\}$ and replacing $c$ with $-c$,
eq.~\eqref{eq:StandardLP.strongDuality:minmin} can be equivalently transformed to
\begin{equation}
\max \,\{\, c \* x ~|~ x \ge 0 \aand A \* x \le b \,\}
= \min \,\{\, b \* y ~|~ y \ge 0 \aand A^T \!\* y \ge c \,\}
\label{eq:StandardLP.strongDuality:maxmin}
\end{equation}
which is probably a more familiar formulation of strong LP duality.
Later it will be clear why we chose the ``$\min\,$/$\,\min$''
rather than the more idiomatic ``$\max\,$/$\,\min$'' formulation. Also, we do not investigate
``asymmetric versions'' such as:
\begin{equation}
\max \,\{\, c \* x ~|~ x \ge 0 \aand A \* x = b \,\}
= \min \,\{\, b \* y ~|~ A^T \!\* y \ge c \,\}
\label{eq:StandardLP.strongDuality:assym}
\end{equation}

For many more Farkas-like theorems (beyond the scope of this paper),
see for example \cite{Perng2017}.

\vspace*{0.001mm}

\noindent This paper makes the following contributions.
\begin{itemize}
\item We formally prove several existing duality results (including the three theorems above) in Lean 4.
In fact, we prove a more general version of \ref{equalityFarkas} due to Bartl~\cite{Bartl2007}.
\sekt{generalizations} says more about it.
\item We establish (and formally prove in Lean 4) a new generalization
of~\ref{inequalityFarkas} and~\ref{StandardLP.strongDuality}
to the case when some of the coefficients are allowed to have infinite values. 
This scenario can be motivated by discrete optimization problems with ``hard'' constraints
(the ``hard'' constraints declare what has to be be satisfied,
whereas ``soft'' constraints declare what should be optimized).
A common way to concisely write down such problems mathematically is to use
infinite coefficients in front of the corresponding terms in the objective. 
Since infinities are not allowed in traditional LPs,
``soft'' and ``hard'' constraints need to be handled differently 
when formulating LP relaxations of such problems (see e.g.~\cite{BLP}).
Our work provides a more direct way to formulate such relaxations.
\sekt{extensions} says more about the extensions.
\end{itemize}

\subsection{Bartl's generalization}\label{generalizations}

The next theorem generalizes \ref{equalityFarkas} to structures where
multiplication does not have to be commutative.
Furthermore, it supports infinitely many equations.

\begin{theorem}[\customlabel{coordinateFarkas}{\texttt{coordinateFarkasBartl}}]
Let $I$ be any type.
Let $J$ be a finite type.
Let $R$ be a linearly ordered division ring.
Let $A$ be an $R$-linear map from $(I \r R)$ to $(J \r R)$.
Let $b$ be an $R$-linear map from $(I \r R)$ to $R$.
Exactly one of the following exists:
\begin{itemize}
\item nonnegative vector $x : J \r R$ such that, for all $w : I \r R$, we have
$ \sum_{j : J}\; (A~w)_j \cdot x_j = b~w $
\item vector $y : I \r R$ such that $A~y \ge 0$ and $b~y < 0$
\end{itemize}
\end{theorem}

Note that \ref{equalityFarkas} for matrix $A: (I \times J) \r F$ and vector $b : I \r F$
can be obtained by applying \ref{coordinateFarkas} to the $F$-linear maps
$(A^T \* \phantom{v})$ and $(b \* \phantom{v})$ utilizing the fact that
two linear maps are equal if and only if they map the basis vectors equally.

In the next generalization (a similar theorem is in \cite{Chernikov}), 
the partially ordered module $I \r R$ is replaced by a general $R$-module $W$.

\begin{theorem}[\customlabel{scalarFarkas}{\texttt{almostFarkasBartl}}]
Let $J$ be a finite type.
Let $R$ be a linearly ordered division ring.
Let $W$ be an $R$-module.
Let $A$ be an $R$-linear map from $W$ to $(J \r R)$.
Let $b$ be an $R$-linear map from $W$ to $R$.
Exactly one of the following exists:
\begin{itemize}
\item nonnegative vector $x : J \r R$ such that, for all $w : W$, we have
$ \sum_{j : J}\; (A~w)_j \cdot x_j = b~w $
\item vector $y : W$ such that $A~y \ge 0$ and $b~y < 0$
\end{itemize}
\end{theorem}
In the most general theorem, stated below, certain occurrences of $R$ are replaced by
a linearly ordered $R$-module $V$ whose order respects the order on $R$
(for a formal definition, see the end of \sekt{preliminaries-modules}).

\begin{theorem}[\customlabel{fintypeFarkasBartl}{\texttt{fintypeFarkasBartl}}]
Let $J$ be a finite type.
Let $R$ be a linearly ordered division ring.
Let $W$ be an $R$-module.
Let $V$ be a linearly ordered $R$-module whose ordering satisfies
monotonicity of scalar multiplication by nonnegative elements on the left.
Let $A$ be an $R$-linear map from $W$ to $(J \r R)$.
Let $b$ be an $R$-linear map from $W$ to $V$.
Exactly one of the following exists:
\begin{itemize}
\item nonnegative vector family $x : J \r V$ such that, for all $w : W$, we have
$ \sum_{j : J}\; (A~w)_j \cdot x_j = b~w $
\item vector $y : W$ such that $A~y \ge 0$ and $b~y < 0$
\end{itemize}
\end{theorem}
In the last branch, $A~y \ge 0$ uses the partial order on $(J \r R)$ whereas
$b~y < 0$ uses the linear order on $V$.
Note that \ref{fintypeFarkasBartl} subsumes \ref{scalarFarkas} (as well as the other versions based on equality),
since $R$ can be viewed as a linearly ordered module over itself.
We prove \ref{fintypeFarkasBartl} in \sekt{bartl}, which is where the heavy lifting comes.
Our proof is based on \cite{Bartl2011}.

\subsection{Extension to infinite coefficients}\label{extensions}

Until now, we have talked about known results.
What follows is a new extension of the theory.

\begin{definition}
Let $F$ be a linearly ordered field.
We define an \emph{extended} linearly ordered field $F_\infty$ as
$F \cup \{ \bot, \top \}$ with the following properties.
Let $p$ and $q$ be numbers from $F$.
We have $\bot < p < \top$.
We define addition, negation, and scalar action on $F_\infty$ as follows:
\begin{center}
	\begin{tabular}{ | c " c | c | c | }
		\hline
		$+$ & $\bot$ & $q$ & $\top$  \\
		\thickhline
		$\bot$ & $\bot$ & $\bot$ & $\bot$  \\ 
		\hline
		$p$ & $\bot$ & $p\!+\!q$ & $\top$  \\ 
		\hline
		$\top$ & $\bot$ & $\top$ & $\top$ \\ 
		\hline
	\end{tabular}
	\qquad\qquad\qquad
	\begin{tabular}{ | c " c | c | c | }
		\hline
		$-$ & $\bot$ & $q$ & $\top$  \\
		\thickhline
		$=$ & $\top$ & $-q$ & $\bot$ \\
		\hline
	\end{tabular}
	\qquad\qquad\qquad
	\begin{tabular}{ | c " c | c | c | }
		\hline
		$\cdot$ & $\bot$ & $q$ & $\top$  \\
		\thickhline
		$0$ & $\bot$ & $0$ & $0$  \\ 
		\hline
		$p>0$ & $\bot$ & $p \cdot q$ & $\top$  \\ 
		\hline
	\end{tabular}
\end{center}
When we talk about elements of $F_\infty$,
we say that values from $F$ are \emph{finite}.
\end{definition}

Informally speaking, $\top$ represents the positive infinity,
$\bot$ represents the negative infinity, and we say that
$\bot$ is ``stronger'' than $\top$ in~all arithmetic operations.
The surprising parts are $\bot + \top = \bot$ and $0 \.\cdot \bot = \bot$.
Because of them, $F_\infty$ is not a field.
In fact, $F_\infty$ is not even a group. 
However, $F_\infty$ is still a densely linearly ordered abelian monoid
with characteristic zero.
Note that $\bot + \top = \bot$ is a standard convention in Mathlib \cite{Mathlib}
but $0 \.\cdot \bot = \bot$ is ad hoc.

\begin{theorem}[\customlabel{extendedFarkas}{\texttt{extendedFarkas}}]
Let $I$ and $J$ be finite types.
Let $F$ be a linearly ordered field.
Let $A$ be a matrix of type $(I \times J) \r F_\infty$.
Let $b$ be a vector of type $I \r F_\infty$.
Assume that $A$ does not have $\bot$ and $\top$ in the same row.
Assume that $A$ does not have $\bot$ and $\top$ in the same column.
Assume that $A$ does not have $\top$ in any row where $b$ has $\top$.
Assume that $A$ does not have $\bot$ in any row where $b$ has~$\bot$.
Exactly one of the following exists:
\begin{itemize}
\item nonnegative vector $x : J \r F$ such that $A \* x \le b$
\item nonnegative vector $y : I \r F$ such that $(-A^T) \* y \le 0$ and $b \* y < 0$
\end{itemize}
\end{theorem}
Note \ref{extendedFarkas} has four preconditions on matrix $A$ and vector $b$. 
In \sekt{sec:extendedFarkas:counterexamples}
we show that omitting any of them makes the theorem false.
Observe that \ref{inequalityFarkas} has condition $A^T \* y \ge 0$
in the second branch, while in \ref{extendedFarkas}
we changed it to  $(-A^T) \* y \le 0$.
The two conditions are equivalent
for finite-valued matrices $A$,
but not necessarily for matrices $A$ with infinities
(e.g.\ condition
$(-\top) \* 0 \ge 0$ is false but 
$\top \* 0 \le 0$ is true).
One must be careful when formulating this condition;
for example, using the condition from \ref{inequalityFarkas}
would make the theorem false even if $A$ has only a single $\bot$ entry
(see \sekt{sec:extendedFarkas:counterexamples}).

Next, we define an extended notion of linear program, i.e.,
linear programming over extended linearly ordered fields.
The implicit intention is that the linear program is to be minimized.

\begin{definition}
Let $I$ and $J$ be finite types.
Let $F$ be a linearly ordered field.
Let $A$ be a matrix of type $(I \times J) \r F_\infty$,
let $b$ be a vector of type $I \r F_\infty$,
and $c$ be a vector of type $J \r F_\infty$.
We say that $P = (A, b, c)$ is a \emph{linear program} over $F_\infty$
whose constraints are indexed by $I$ and variables are indexed by $J$.

A nonnegative vector $x : J \r R$ is a \emph{solution} to $P$ if  $A \* x \le b$.
We say that $P$ is \emph{feasible} if there exists solution $x$ with $c \* x \ne \top$.
We say that $P$ is \emph{unbounded} if, for any $r\in F$, there exists solution $x$
with $c \* x \le r$.

The \emph{optimum} of $P$, denoted as $P^\star$,
is defined as follows. If $P$ is not feasible then $P^\star=\top$.
Else, if $P$ is unbounded, then $P^\star=\bot$.
Else, if there exists finite $r\in F$ such that $c \* x \ge r$ for all solutions $x$
and $c \* x^\star = r$ for some solution $x^\star$, then $P^\star=r$.
Otherwise, $P^\star$ is undefined.
\end{definition}

In order to state our duality results, we need a few more definitions.

\begin{definition}
Linear Program $P = (A, b, c)$ over $F_\infty$ is said to be
{\em valid} if it satisfies the following six conditions:
\begin{itemize}
\item $A$ does not have $\bot$ and $\top$ in the same row
\item $A$ does not have $\bot$ and $\top$ in the same column
\item $A$ does not have $\bot$ in any row where $b$ has $\bot$
\item $A$ does not have $\top$ in any column where $c$ has $\bot$
\item $A$ does not have $\top$ in any row where $b$ has $\top$
\item $A$ does not have $\bot$ in any column where $c$ has $\top$
\end{itemize}
We say that the linear program $(-A^T, c, b)$ is the \emph{dual} of $P$.
\end{definition}
It is straightforward to check that the dual of a valid LP is valid.
Later we show that a valid $P$ cannot have $P^\star$ undefined.

\begin{theorem}[\customlabel{ExtendedLP.strongDuality}{\texttt{ValidELP.strongDuality}}]
Let $F$ be a linearly ordered field,
$P$ be a valid linear program over $F_\infty$ and $D$ be its dual.
Then $P^\star$ and $D^\star$ are defined.
If at least one of them is feasible, i.e., $(P^\star, D^\star) \ne (\top, \top)$,
then $P^\star = -D^\star$.
\end{theorem}
\begin{sloppypar} 
Similar to~\ref{extendedFarkas}, all six conditions in the definition of a valid LP are necessary;
omitting any one of them makes \ref{ExtendedLP.strongDuality} false
(see counterexamples in \sekt{sec:ExtendedLP.strongDuality:counterexamples}).
\end{sloppypar}

Note that the conclusion of the theorem ($P^\star = -D^\star$) can be reformulated,
when we use highly informal notation, as in eq.~\eqref{eq:StandardLP.strongDuality:minmin}:
$$ \min \,\{\, c \* x ~|~ x \ge 0 \aand A \* x \le b \,\}
=- \min \,\{\, b \* y ~|~ y \ge 0 \aand (-A^T) \* y \le c \,\}
$$
If all entries of $(A,b,c)$ are finite then this equation can be stated in
many equivalent ways, e.g.\ as in~\eqref{eq:StandardLP.strongDuality:maxmin}.
This reformulation uses the facts
that $ (-c) \* x = - (c \* x) $,
and that $(-A^T) \* y \le c$ is equivalent to $A^T \* y \ge -c$.
These facts are no longer true if infinities are allowed,
so one must  be careful when formulating the duality theorem for $F_\infty$.
We considered several ``$\max\,$/$\,\min$'' and ``$\max\,$/$\,\max$'' versions,
but they all required stronger preconditions compared to the ``$\min\,$/$\,\min$'' version
given in~\ref{ExtendedLP.strongDuality}.

\subsubsection{Example---cheap lunch}

According to Nutritionix\footnote{\url{https://www.nutritionix.com/food/white-rice/1000-g}} a kilogram of boiled white rice contains $27$ g protein and $1300$ kcal.
According to Nutritionix\footnote{\url{https://www.nutritionix.com/food/lentils}} a kilogram of boiled lentils contains $90$ g protein and $1150$ kcal.
Let's say that a kilogram of boiled white rice costs $0.92$ euro and that
a kilogram of boiled lentils costs $1.75$ euro.
We want to cook a lunch, as cheap as possible, that contains at least
$30$ g protein and at least $700$ kcal.
The choice of white rice and lentils isn't random\:---\:one of the authors
ate this lunch at a mathematical camp\:---\:at the moment, he didn't like the
lunch\:---\:it wasn't very tasty\:---\:but later he realized what an awesome nutritional
value the lunch had for its low price.
This is a simple example of the well-known diet problem \cite{Dantzig}.
$$
\begin{array}{rcl}
\min && 0.92 \cdot r + 1.75 \cdot l \\
(-27) \cdot r + (-90) \cdot l &\le& - 30 \\
(-1300) \cdot r + (-1150) \cdot l &\le& - 700 \\
r &\ge& 0  \\
l &\ge& 0 
\end{array}
$$
Using a numerical LP solver, we obtain the optimal solution
$r=0.331588$ and $l=0.233857$ which satisfies our dietary requirements for
$0.714311$ euro.
The dual of this problem, in the sense of \ref{ExtendedLP.strongDuality}, is the following LP:
$$
\begin{array}{rcl}
\min && (-30) \cdot p + (-700) \cdot k \\
27 \cdot p + 1300 \cdot k &\le& 0.92 \\
 90 \cdot p + 1150 \cdot k &\le& 1.75 \\
p &\ge& 0  \\
k &\ge& 0 
\end{array}
$$
Using a numerical LP solver, we obtain the optimal solution
$p=0.0141594$ and $k=0.000413613$ leading to the
objective value $-0.714311$ here.
To summarize, the cheapest lunch that contains at least
30 g protein and at least 700 kcal consists of $332$ g white rice
and $234$ g lentils; the shadow cost of a gram of protein is
$0.014$ euro and the shadow cost of a kcal is $0.00041$ euro
in our settings.

Now consider settings where lentils are out of stock.
We still want to obtain at least 30 g protein and at least 700 kcal.
What is the price of our lunch now?
On paper, it is most natural to entirely remove lentils out of the picture
and recompute the LP.
However, computer proof systems generally don't like it when
sizes of matrices change, so let's do the modification of input in place.
One could suggest that, in order to model the newly arised situation,
the price of lentils will be increased to $999999$ euro,
so that the optimal solution will assign $0$ to it.
This is, however, not mathematically elegant, as it requires some
engineering insight into setting the constant to an appropriately large number.
We claim that the proper mathematical approach is to increase the
price of lentils to infinity! Let's see how our original LP will look now:
$$
\begin{array}{rcl}
\min && 0.92 \cdot r + \top \cdot l \\
(-27) \cdot r + (-90) \cdot l &\le& - 30 \\
(-1300) \cdot r + (-1150) \cdot l &\le& - 700 \\
r &\ge& 0  \\
l &\ge& 0 
\end{array}
$$
The optimal solution is now $r=1.111111$ and $l=0$ and costs $1.022222$ euro.
The dual of this problem is the following LP:
$$
\begin{array}{rcl}
\min && (-30) \cdot p + (-700) \cdot k \\
27 \cdot p + 1300 \cdot k &\le& 0.92 \\
 90 \cdot p + 1150 \cdot k &\le& \top \\
p &\ge& 0  \\
k &\ge& 0 
\end{array}
$$
At this point, the shadow price of a gram of protein is $0.034074$ euro but the shadow price of kcal is $0$ euro.
The objective value is $-1.022222$ in accordance with our theoretical prediction.

\subsection{Structure of this paper}

In \sekt{preliminaries}, we explain all underlying definitions and
comment on the formalization process;
following the philosophy of \cite{Believe}
we review most of the declarations needed for the reader to believe
our results, leaving out many declarations that were used in order
to prove the results.
In \sekt{statements}, we formally state theorems from \sekt{introduction}
using definitions from \sekt{preliminaries}.
In \sekt{bartl}, we prove \ref{fintypeFarkasBartl} (stated in \sekt{generalizations}),
from which we obtain \ref{equalityFarkas} as a corollary.
In \sekt{extended}, we prove \ref{extendedFarkas} (stated in \sekt{extensions}).
In \sekt{dualityELP}, we prove \ref{ExtendedLP.strongDuality}
(stated in \sekt{extensions}), from which we obtain
\ref{StandardLP.strongDuality} as a corollary.
In \sekt{counterexamples}, we show what happens when various preconditions
are not satisfied.

Repository \texttt{https://github.com/madvorak/duality/tree/v3.2}
contains the full version of all definitions, statements,
and proofs. They are written in a formal language called
Lean, which provides a guarantee that every step of
every proof follows from valid logical axioms.\footnote{
The only axioms used in our proofs are
\texttt{propext}, \texttt{Classical.choice}, and
\texttt{Quot.sound}, which you can check by the
\texttt{\#print axioms} command.}
We use Lean version 4.18.0 together with the Lean
mathematical library \texttt{Mathlib} \cite{Mathlib}
revision \texttt{aa936c3} (dated 2025-04-01).
This paper attempts to be an accurate description of
the \texttt{duality} project.
However, in case of any discrepancy, the code shall prevail.

\section{Preliminaries}
\label{preliminaries}

There are many layers of definitions that are built before say
what a linearly ordered field is,
what a linearly ordered division ring is, and
what a linearly ordered abelian group (used in \ref{fintypeFarkasBartl}) is.
\sekt{preliminaries-typeclasses} mainly documents existing Mathlib definitions
(but omitting declarations of default values);
we also highlight how we added linearly ordered division rings into our project.
\sekt{preliminaries-ef} then explains how we defined extended linearly ordered fields.
\sekt{preliminaries-vectors} reviews vectors and matrices; we especially focus on
how multiplication between them is defined.
\sekt{preliminaries-modules} explains how modules are implemented.
\sekt{preliminaries-lp} provides the formal definition of extended linear programs.

\subsection{Review of algebraic typeclasses that our project depends on}
\label{preliminaries-typeclasses}

Additive semigroup is a structure on any type with addition (denoted by the infix \texttt{+} operator) where the addition is associative:
\begin{lstlisting}
class AddSemigroup (G : Type u) extends Add G where
  add_assoc : ∀ a b c : G, (a + b) + c = a + (b + c)
\end{lstlisting}
Semigroup is a structure on any type with multiplication (denoted by the infix \texttt{*} operator) where the multiplication is associative:
\begin{lstlisting}
class Semigroup (G : Type u) extends Mul G where
  mul_assoc : ∀ a b c : G, (a * b) * c = a * (b * c)
\end{lstlisting}
Additive monoid is an additive semigroup with the ``zero'' element that is neutral with respect to addition
from both left and right, equipped with a scalar multiplication by the natural numbers:
\begin{lstlisting}
class AddZeroClass (M : Type u) extends Zero M, Add M where
  zero_add : ∀ a : M, 0 + a = a
  add_zero : ∀ a : M, a + 0 = a
class AddMonoid (M : Type u) extends AddSemigroup M, AddZeroClass M where
  nsmul : ℕ → M → M
  nsmul_zero : ∀ x : M, nsmul 0 x = 0 
  nsmul_succ : ∀ (n : ℕ) (x : M), nsmul (n + 1) x = nsmul n x + x 
\end{lstlisting}
Similarly, monoid is a semigroup with the ``one'' element that is neutral with respect to multiplication
from both left and right, equipped with a power to the natural numbers:
\begin{lstlisting}
class MulOneClass (M : Type u) extends One M, Mul M where
  one_mul : ∀ a : M, 1 * a = a
  mul_one : ∀ a : M, a * 1 = a
class Monoid (M : Type u) extends Semigroup M, MulOneClass M where
  npow : ℕ → M → M
  npow_zero : ∀ x : M, npow 0 x = 1 
  npow_succ : ∀ (n : ℕ) (x : M), npow (n + 1) x = npow n x * x 
\end{lstlisting}
Subtractive monoid is an additive monoid that adds two more operations (unary and binary minus)
that satisfy some basic properties (please note that ``adding minus itself gives zero''
is not required yet; that will be required e.g.~in an additive group):
\begin{lstlisting}
class SubNegMonoid (G : Type u) extends AddMonoid G, Neg G, Sub G where
  sub_eq_add_neg : ∀ a b : G, a - b = a + -b 
  zsmul : ℤ → G → G
  zsmul_zero' : ∀ a : G, zsmul 0 a = 0 
  zsmul_succ' (n : ℕ) (a : G) : zsmul n.succ a = zsmul n a + a
  zsmul_neg' (n : ℕ) (a : G) : zsmul (Int.negSucc n) a = -(zsmul n.succ a)
\end{lstlisting}
Similarly, division monoid is a monoid that adds two more operations (inverse and divide)
that satisfy some basic properties (please note that ``multiplication by an inverse gives one''
is not required yet):
\begin{lstlisting}
class DivInvMonoid (G : Type u) extends Monoid G, Inv G, Div G where
  div_eq_mul_inv : ∀ a b : G, a / b = a * b⁻¹ 
  zpow : ℤ → G → G
  zpow_zero' : ∀ a : G, zpow 0 a = 1 
  zpow_succ' (n : ℕ) (a : G) : zpow n.succ a = zpow n a * a
  zpow_neg' (n : ℕ) (a : G) : zpow (Int.negSucc n) a = (zpow n.succ a)⁻¹ 
\end{lstlisting}
Additive group is a subtractive monoid in which the unary minus acts as a left inverse with respect to addition:
\begin{lstlisting}
class AddGroup (A : Type u) extends SubNegMonoid A where
  neg_add_cancel : ∀ a : A, -a + a = 0
\end{lstlisting}
Abelian magma is a structure on any type that has commutative addition:
\begin{lstlisting}
class AddCommMagma (G : Type u) extends Add G where
  add_comm : ∀ a b : G, a + b = b + a
\end{lstlisting}
Similarly, commutative magma is a structure on any type that has commutative multiplication:
\begin{lstlisting}
class CommMagma (G : Type u) extends Mul G where
  mul_comm : ∀ a b : G, a * b = b * a
\end{lstlisting}
Abelian semigroup is an abelian magma and an additive semigroup at the same time:
\begin{lstlisting}
class AddCommSemigroup (G : Type u) extends AddSemigroup G, AddCommMagma G
\end{lstlisting}
Similarly, commutative semigroup is a commutative magma and a semigroup at the same time:
\begin{lstlisting}
class CommSemigroup (G : Type u) extends Semigroup G, CommMagma G
\end{lstlisting} \pagebreak[2]
Abelian monoid is an additive monoid and an abelian semigroup at the same time:
\begin{lstlisting}
class AddCommMonoid (M : Type u) extends AddMonoid M, AddCommSemigroup M
\end{lstlisting}
Similarly, commutative monoid is a monoid and a commutative semigroup at the same time:
\begin{lstlisting}
class CommMonoid (M : Type u) extends Monoid M, CommSemigroup M
\end{lstlisting}
Abelian group is an additive group and an abelian monoid at the same time:
\begin{lstlisting}
class AddCommGroup (G : Type u) extends AddGroup G, AddCommMonoid G
\end{lstlisting}
Distrib is a structure on any type with addition and multiplication where
both left distributivity and right distributivity hold
(please ignore the difference between \texttt{Type u} and \texttt{Type*} as it is the same thing
for all our purposes):
\begin{lstlisting}
class Distrib (R : Type*) extends Mul R, Add R where
  left_distrib : ∀ a b c : R, a * (b + c) = a * b + a * c
  right_distrib : ∀ a b c : R, (a + b) * c = a * c + b * c
\end{lstlisting}
Nonunital-nonassociative-semiring is an abelian monoid with distributive multiplication and well-behaved zero:
\begin{lstlisting}
class MulZeroClass (M₀ : Type u) extends Mul M₀, Zero M₀ where
  zero_mul : ∀ a : M₀, 0 * a = 0
  mul_zero : ∀ a : M₀, a * 0 = 0
class NonUnitalNonAssocSemiring (α : Type u) extends AddCommMonoid α, Distrib α, MulZeroClass α
\end{lstlisting}
Nonunital-semiring is a nonunital-nonassociative-semiring that forms a semigroup with zero (i.e., the
semigroup-with-zero requirement makes it associative):
\begin{lstlisting}
class SemigroupWithZero (S₀ : Type u) extends Semigroup S₀, MulZeroClass S₀
class NonUnitalSemiring (α : Type u) extends NonUnitalNonAssocSemiring α, SemigroupWithZero α
\end{lstlisting}
Additive monoid with one and abelian monoid with one are defined from additive monoid
and abelian monoid, equipped with the symbol ``one'' and embedding of natural numbers:
\begin{lstlisting}
class AddMonoidWithOne (R : Type*) extends NatCast R, AddMonoid R, One R where
  natCast_zero : natCast 0 = 0
  natCast_succ : ∀ n : ℕ, natCast (n + 1) = (natCast n) + 1 
class AddCommMonoidWithOne (R : Type*) extends AddMonoidWithOne R, AddCommMonoid R
\end{lstlisting}
Additive group with one is an additive monoid with one and additive group, and embeds all integers:
\begin{lstlisting}
class AddGroupWithOne (R : Type u) extends IntCast R, AddMonoidWithOne R, AddGroup R where
  intCast_ofNat : ∀ n : ℕ, intCast (n : ℕ) = Nat.cast n 
  intCast_negSucc : ∀ n : ℕ, intCast (Int.negSucc n) = - Nat.cast (n + 1) 
\end{lstlisting}
Nonassociative-semiring is a nonunital-nonassociative-semiring that has a well-behaved multiplication
by both zero and one and forms an abelian monoid with one (i.e., the unit is finally defined):
\begin{lstlisting}
class MulZeroOneClass (M₀ : Type u) extends MulOneClass M₀, MulZeroClass M₀
class NonAssocSemiring (α : Type u) extends NonUnitalNonAssocSemiring α, MulZeroOneClass α,
  AddCommMonoidWithOne α
\end{lstlisting}
Semiring is a nonunital-semiring and nonassociative-semiring that forms a monoid with zero:
\begin{lstlisting}
class MonoidWithZero (M₀ : Type u) extends Monoid M₀, MulZeroOneClass M₀, SemigroupWithZero M₀
class Semiring (α : Type u) extends NonUnitalSemiring α, NonAssocSemiring α, MonoidWithZero α
\end{lstlisting}
Ring is a semiring and an abelian group at the same time that has ``one'' that behaves well:
\begin{lstlisting}
class Ring (R : Type u) extends Semiring R, AddCommGroup R, AddGroupWithOne R
\end{lstlisting}
Commutative ring is a ring (guarantees commutative addition) and a commutative monoid
(guarantees commutative multiplication) at the same time:
\begin{lstlisting}
class CommRing (α : Type u) extends Ring α, CommMonoid α
\end{lstlisting}
Division ring is a nontrivial ring whose multiplication forms a division monoid, whose nonzero elements have
multiplicative inverses, whose zero is inverse to itself
(if you find the equality $0^{-1}\!=\!0$ disturbing, read \cite{DivBy0} that explains it),
and embeds rational numbers:
\begin{lstlisting}
class Nontrivial (α : Type*) : Prop where
  exists_pair_ne : ∃ x y : α, x ≠ y
class DivisionRing (K : Type*) extends Ring K, DivInvMonoid K, Nontrivial K, NNRatCast K, RatCast K where
  mul_inv_cancel : ∀ (a : K), a ≠ 0 → a * a⁻¹ = 1
  inv_zero : (0 : K)⁻¹ = 0
\end{lstlisting}
Field is a commutative ring and a division ring at the same time:
\begin{lstlisting}
class Field (K : Type u) extends CommRing K, DivisionRing K
\end{lstlisting}
Preorder is a reflexive \& transitive relation on any structure with binary relational symbols $\le$ and $<$
where the strict comparison $a < b$ is equivalent to $a \le b \aand \neg (b \le a)$ given by the relation $\le$
which is neither required to be symmetric nor required to be antisymmetric:
\begin{lstlisting}
class Preorder (α : Type u) extends LE α, LT α where
  le_refl : ∀ a : α, a ≤ a
  le_trans : ∀ a b c : α, a ≤ b → b ≤ c → a ≤ c
  lt_iff_le_not_le : ∀ a b : α, a < b  ↔  a ≤ b ∧ ¬b ≤ a 
\end{lstlisting}
Partial order is an antisymmetric preoder (that is, a reflexive \& antisymmetric \& transitive relation):
\begin{lstlisting}
class PartialOrder (α : Type u) extends Preorder α where
  le_antisymm : ∀ a b : α, a ≤ b → b ≤ a → a = b
\end{lstlisting}
Ordered abelian group is an abelian group with partial order that respects addition:
\begin{lstlisting}
class OrderedAddCommGroup (α : Type u) extends AddCommGroup α, PartialOrder α where
  add_le_add_left : ∀ a b : α, a ≤ b → ∀ c : α, c + a ≤ c + b
\end{lstlisting}
Strictly ordered ring is a nontrivial ring whose addition behaves as an ordered abelian group
where zero is less or equal to one and the product of two strictly positive elements is strictly positive:
\begin{lstlisting}
class StrictOrderedRing (α : Type u) extends Ring α, OrderedAddCommGroup α, Nontrivial α where
  zero_le_one : 0 ≤ (1 : α)
  mul_pos : ∀ a b : α, 0 < a → 0 < b → 0 < a * b
\end{lstlisting}
Linear order (sometimes called total order) is a partial order where every two elements are comparable;
technical details are omitted:
\begin{lstlisting}
class LinearOrder (α : Type u) extends PartialOrder α, (...)
  le_total (a b : α) : a ≤ b ∨ b ≤ a
  (...)
\end{lstlisting}
Linearly ordered abelian group is an ordered abelian group whose order is linear:
\begin{lstlisting}
class LinearOrderedAddCommGroup (α : Type u) extends OrderedAddCommGroup α, LinearOrder α
\end{lstlisting}
Linearly ordered ring is a strictly ordered ring where every two elements are comparable:
\begin{lstlisting}
class LinearOrderedRing (α : Type u) extends StrictOrderedRing α, LinearOrder α
\end{lstlisting}
Linearly ordered commutative ring is a linearly ordered ring and commutative monoid at the same time:
\begin{lstlisting}
class LinearOrderedCommRing (α : Type u) extends LinearOrderedRing α, CommMonoid α
\end{lstlisting}
In our project, we define a linearly ordered division ring as a division ring that is a linearly ordered ring at the same time:
\begin{lstlisting}
class LinearOrderedDivisionRing (R : Type*) extends LinearOrderedRing R, DivisionRing R
\end{lstlisting}
Linearly ordered field is defined in Mathlib as a linearly ordered commutative ring that is a field at the same time:
\begin{lstlisting}
class LinearOrderedField (α : Type*) extends LinearOrderedCommRing α, Field α
\end{lstlisting}
Note that \texttt{LinearOrderedDivisionRing} is not a part of the algebraic hierarchy provided by Mathlib,
hence \texttt{LinearOrderedField} does~not~inherit \texttt{LinearOrderedDivisionRing}.

Because of that, we provide a custom instance that converts
\texttt{LinearOrderedField} to \texttt{LinearOrderedDivisionRing}~as~follows:
\begin{lstlisting}
instance LinearOrderedField.toLinearOrderedDivisionRing {F : Type*} [instF : LinearOrderedField F] :
  LinearOrderedDivisionRing F := { instF with }
\end{lstlisting}
This instance is needed for the step from \ref{coordinateFarkas} to \ref{equalityFarkas}.

\subsection{Extended linearly ordered fields}
\label{preliminaries-ef}

Given any type $F$, we construct $F \cup \{ \bot, \top \}$ as follows:
\begin{lstlisting}
def Extend (F : Type*) := WithBot (WithTop F)
\end{lstlisting}
From now on we assume that $F$ is a linearly ordered field:
\begin{lstlisting}
variable {F : Type*} [LinearOrderedField F]
\end{lstlisting}
The following instance defines how addition and comparison
behaves on $F_\infty$ and automatically generates a proof
that $F_\infty$ forms a linearly ordered abelian monoid:
\begin{lstlisting}
instance : LinearOrderedAddCommMonoid (Extend F) :=
  inferInstanceAs (LinearOrderedAddCommMonoid (WithBot (WithTop F)))
\end{lstlisting}
The following definition embeds $F$ in $F_\infty$ and registers
this canonical embedding as a type coercion:
\begin{lstlisting}
@[coe] def toE : F → (Extend F) := some ∘ some
instance : Coe F (Extend F) := ⟨toE⟩
\end{lstlisting} \pagebreak[2]
Unary minus on $F_\infty$ is defined as follows:
\begin{lstlisting}
def neg : Extend F → Extend F
| ⊥ => ⊤
| ⊤ => ⊥
| (x : F) => toE (-x)
instance : Neg (Extend F) := ⟨EF.neg⟩
\end{lstlisting}
Line-by-line, we see that:
\begin{itemize}
\item negating $\bot$ gives $\top$
\item negating $\top$ gives $\bot$
\item negating any finite value $x$ gives $-x$ converted to the type $F_\infty$
\end{itemize}
In the file \texttt{FarkasSpecial.lean} and everything downstream,
we have the notation
\begin{lstlisting}
F∞
\end{lstlisting}
for $F_\infty$ and also the notation
\begin{lstlisting}
F≥0
\end{lstlisting}
for the type of nonnegative elements of $F$.
We define a scalar action of the nonnegative elements
on $F_\infty$ as follows:
\begin{lstlisting}
def EF.smulNN (c : F≥0) : F∞ → F∞
| ⊥ => ⊥
| ⊤ => if c = 0 then 0 else ⊤
| (f : F) => toE (c.val * f)
\end{lstlisting}
Line-by-line, we see that:
\begin{itemize}
\item multiplying $\bot$ by anything from the left gives $\bot$
\item multiplying $\top$ by $0$ from the left gives $0$;
multiplying $\top$ by a strictly positive element from the left gives $\top$
\item multiplying any finite value $f$ by a coefficient $c$ from the left
gives $c \* f$ converted to the type $F_\infty$
\end{itemize}
Scalar action on $F_\infty$ by negative elements from the left is undefined.
Multiplication between two elements of $F_\infty$ is also undefined.

\subsection{Vectors and matrices}
\label{preliminaries-vectors}

We distinguish two types of vectors; implicit vectors and explicit vectors.
Implicit vectors (called just ``vectors'') are members of a vector space;
they do not have any internal structure
(in the informal text, we use the word ``vectors'' somewhat loosely;
it can refer to members of any module).
Explicit vectors are functions from coordinates to values (again, we use the word
``explicit vector`` (or just ``vector'' when it is clear from context that our vector
is a map) not only when they form a vector space).
The type of coordinates does not have to be ordered.

Matrices live next to explicit vectors. They are also functions; they take a row index
and a column index and they output a value at the given spot.
Neither the row indices nor the column vertices are required to form an ordered type.
This is why multiplication between matrices and vectors is defined only in structures
where addition forms an abelian monoid. Consider the following example:
$$
\begin{pmatrix}
	1 & 2 & 3 \\
	4 & 5 & 6 \\
\end{pmatrix}
\*
\begin{pmatrix}
	7 \\ 8 \\ 9
\end{pmatrix}
=
\begin{pmatrix}
	? \\ \_
\end{pmatrix}
$$
We do not know whether the value at the question mark is equal to
$ 1 \* 7 + 2 \* 8 + 3 \* 9 $ or to
$ 2 \* 8 + 1 \* 7 + 3 \* 9 $ or to
any other ordering of summands.
This is why commutativity of addition is necessary for the definition to be valid.
On the other hand, we do not assume any property of multiplication in the
definition of multiplication between matrices and vectors; they do not even
have to be of the same type; we only require the elements of the vector
to have an action on the elements of the matrix (this is not a typo\:---\:normally,
we would want matrices to have an action on vectors\:---\:not in our work).
\sekt{versus} provides more details on multiplication between matrices and vectors
and on the dot product.

One thing to be said about the implementation of matrices is that we lied
a little bit when we wrote things like $A : (I \times J) \r F$.
In reality, \texttt{A : Matrix I J F} is implemented as $A : \,I \r (J \r F)$.
This technique \cite{Currying} is called currying and makes it easier to
pass individual rows of the matrix into functions and lemmas.
\pagebreak[2]

The following instance defines how explicit vectors are compared;
Mathlib defines $x \le y$ to hold if and only if
$x_i \le y_i$ holds for~every~$i$ 
(it is, in fact, more general than we said because the following instance
talks about $\Pi$-types \cite{MLTT}, i.e.,
the \textit{type} of \texttt{x i} can depend on the \textit{value} of \texttt{i}):
\newpage
\begin{lstlisting}
instance Pi.hasLe {ι : Type*} {π : ι → Type*} [∀ i, LE (π i)] :
  LE (∀ i, π i) where le x y := ∀ i, x i ≤ y i
\end{lstlisting}
Comparison between matrices is not defined, only equality for matrices is.

Throughout the paper, we do not distinguish between
nonnegative explicit vectors and explicit vectors of nonnegative elements.
The code, however, does distinguish between them.

\subsubsection{Mathlib definitions versus our definitions}
\label{versus}

Mathlib defines dot product (i.e., a product of an explicit vector with an explicit vector) as follows:
\begin{lstlisting}
def dotProduct [Fintype m] [Mul α] [AddCommMonoid α] (v w : m → α) : α :=
  ∑ i, v i * w i
\end{lstlisting}
Mathlib defines a product of a matrix with an explicit vector as follows:
\begin{lstlisting}
def Matrix.mulVec [Fintype n] [NonUnitalNonAssocSemiring α] (M : Matrix m n α) (v : n → α) : m → α
  | i => (fun j => M i j) ⬝ᵥ v
\end{lstlisting}
Mildly confusingly, \texttt{m} and \texttt{n} are type variables here, not natural numbers.
Note that, when multiplying two vectors,
the left vector is not transposed\:---\:a vector is not
defined as a special case of matrix, and transposition is not
defined as an operation on explicit vectors, only on matrices.
Explicit vectors are neither row vectors nor column vectors;
they are just maps. It is only in our imagination that we
treat certain occurrences of vectors as rows or columns;
in reality, it is the position of the vector in the term that
determines what the vector does there.

Infix notation for these two operations is defined as follows
(the keyword \texttt{infixl} when defining an operator $\circ$ means that
\hbox{\texttt{x $\circ$ y $\circ$ z}} gets parsed as \texttt{(x $\circ$ y) $\circ$ z}
whether or not it makes sense, whereas the keyword \texttt{infixr} when defining an operator $\circ$ means that 
\hbox{\texttt{A $\circ$ B $\circ$ v}} gets parsed as \texttt{A $\circ$ (B $\circ$ v)}
whether or not it makes sense; the numbers 72 and 73 determine precedence\:---\:higher number means tighter):
\begin{lstlisting}
infixl:72 " ⬝ᵥ " => dotProduct
scoped infixr:73 " *ᵥ " => Matrix.mulVec
\end{lstlisting}
The definitions above are sufficient for stating results based on linearly ordered fields.
However, our results from \sekt{extensions} require a more general notion of
multiplication between matrices and vectors.
We define them in the section \texttt{hetero\_matrix\_products\_defs}
as follows:
\newpage
\begin{lstlisting}
variable {α γ : Type*} [AddCommMonoid α] [SMul γ α]

def dotWeig (v : I → α) (w : I → γ) : α :=
  ∑ i : I, w i • v i
infixl:72 " ᵥ⬝ " => dotWeig

def Matrix.mulWeig (M : Matrix I J α) (w : J → γ) (i : I) : α :=
  M i ᵥ⬝ w
infixr:73 " ₘ* " => Matrix.mulWeig
\end{lstlisting}
We start by declaring that $\alpha$ and $\gamma$ are types from
any universe (not necessarily both from the same universe).
We require that $\alpha$ forms an abelian monoid and that
$\gamma$ has a scalar action on $\alpha$. In this setting,
we can instantiate $\alpha$ with $F_\infty$ and $\gamma$
with $F_{\ge 0}$ for any linearly ordered field $F$.

For explicit vectors
$v : I \r \alpha$ and $w : I \r \gamma$, we define
their product of type $\alpha$ as follows.
Every element of $v$ gets multiplied from left by
an element of $w$ on the same index.
Then we sum them all together (in unspecified order).
For a matrix $M : (I \times J) \r \alpha$ and
a vector $w : J \r \gamma$, we define
their product of type $I \r \alpha$ as a function
that takes an index $i$ and outputs the dot product
between the $i$-th row of $M$ and the vector $w$.

Beware that the arguments (both in the function definition and
in the infix notation) come in the opposite order from how
scalar action is written. We recommend a mnemonic 
``vector times weights'' for $v \* w$ and
``matrix times weights'' for $M \* w$ where
arguments come in alphabetical order.

In the infix notation, you can distinguish between the standard
Mathlib definitions and our definitions by observing that Mathlib
operators put the letter $_v$ to the right of the symbol whereas
our operators put a letter to the left of the symbol.

Since we have new definitions, we have to rebuild all API (a lot of lemmas)
for \texttt{dotWeig} and \texttt{Matrix.mulWeig} from scratch.
This process was very tiresome. We decided not to develop a full reusable
library, but prove only those lemmas we wanted to use in our project.
For similar reasons, we did not generalize the Mathlib definition of
``vector times matrix'', as ``matrix times vector'' was all we needed.
It was still a lot of lemmas.

\subsection{Modules and how to order them}
\label{preliminaries-modules}

Given types $\alpha$ and $\beta$ such that $\alpha$ has a scalar action on $\beta$
(denoted by the infix • operator)
and $\alpha$ forms a monoid, Mathlib defines multiplicative action where
$1$ of type $\alpha$ gives the identity action and multiplication in the monoid
associates with the scalar action:
\newpage
\begin{lstlisting}
class MulAction (α : Type*) (β : Type*) [Monoid α] extends SMul α β where
  one_smul : ∀ b : β, (1 : α) • b = b
  mul_smul : ∀ (x y : α) (b : β), (x * y) • b = x • y • b
\end{lstlisting}
For a distributive multiplicative action, we furthermore require the latter type
to form an additive monoid and two more properties are required; multiplying the zero
gives the zero, and the multiplicative action is distributive with respect to the addition:
\begin{lstlisting}
class DistribMulAction (M A : Type*) [Monoid M] [AddMonoid A] extends MulAction M A where
  smul_zero : ∀ a : M, a • (0 : A) = 0
  smul_add : ∀ (a : M) (x y : A), a • (x + y) = a • x + a • y
\end{lstlisting}
We can finally review the definition of a module. Here the former type must form a semiring
and the latter type abelian monoid. Module requires a distributive multiplicative action and
two additional properties; addition in the semiring distributes with the multiplicative action,
and multiplying by the zero from the semiring gives a zero in the abelian monoid:
\begin{lstlisting}
class Module (R : Type u) (M : Type v) [Semiring R] [AddCommMonoid M] extends DistribMulAction R M where
  add_smul : ∀ (r s : R) (x : M), (r + s) • x = r • x + s • x
  zero_smul : ∀ x : M, (0 : R) • x = 0
\end{lstlisting}
Note the class \texttt{Module} does not extend the class \texttt{Semiring}; instead,
it requires \texttt{Semiring} as an argument.
The abelian monoid is also required as argument in the definition. We call such a class ``mixin''.
Thanks to this design, we do not need to define subclasses of \texttt{Module} in order to require
``more than a module''. Instead, we use subclasses in the respective arguments, i.e.,
``more than a semiring'' and/or ``more than an abelian monoid''. For example, if we replace
\begin{lstlisting}
[Semiring R] [AddCommMonoid M] [Module R M]
\end{lstlisting}
in a theorem statement by
\begin{lstlisting}
[Field R] [AddCommGroup M] [Module R M]
\end{lstlisting}
we require \texttt{M} to be a vector space over \texttt{R}. We do not need to extend
\texttt{Module} in order to define what a vector space is.

In our case, to state the theorem \ref{fintypeFarkasBartl},
we need \texttt{R} to be a linearly ordered division ring,
we need \texttt{W} to be an \texttt{R}-module, and
we need \texttt{V} to be a linearly ordered \texttt{R}-module.
The following list of requirements is almost correct:
\begin{lstlisting}
[LinearOrderedDivisionRing R] [AddCommGroup W] [Module R W]
[LinearOrderedAddCommGroup V] [Module R V]
\end{lstlisting}
The only missing assumption is the relationship between
how \texttt{R} is ordered and how \texttt{V} is ordered.
For that, we use another mixin, defined in Mathlib as follows:
\begin{lstlisting}
class PosSMulMono (α β : Type*) [SMul α β] [Preorder α] [Preorder β] [Zero α] : Prop where
  elim {a : α} (ha : 0 ≤ a) {b₁ b₂ : β} (hb : b₁ ≤ b₂) : a • b₁ ≤ a • b₂
\end{lstlisting}
Adding \texttt{[PosSMulMono R V]} to the list of requirements solves the issue.

We encourage the reader to try and delete various assumptions and see which parts of the proofs
get underlined in red.

\subsection{Linear programming}
\label{preliminaries-lp}

Extended linear programs are defined by a matrix $A$ and vectors $b$ and $c$
(all of them over an extended linearly ordered field):
\begin{lstlisting}
structure ExtendedLP (I J F : Type*) [LinearOrderedField F] where
  A : Matrix I J F∞
  b : I → F∞
  c : J → F∞
\end{lstlisting}
Vector $x$ made of finite nonnegative values is a solution if and only if $A \* x \le b$
holds:
\begin{lstlisting}
def ExtendedLP.IsSolution [Fintype J] (P : ExtendedLP I J F) (x : J → F≥0) : Prop :=
  P.A ₘ* x ≤ P.b
\end{lstlisting}
$P$ reaches a value $r$ if and only if $P$ has solution $x$ such that $c \* x = r$ holds:
\begin{lstlisting}
def ExtendedLP.Reaches [Fintype J] (P : ExtendedLP I J F) (r : F∞) : Prop :=
  ∃ x : J → F≥0, P.IsSolution x ∧ P.c ᵥ⬝ x = r
\end{lstlisting}
$P$ is feasible if and only if $P$ reaches a value different from $\top$:
\begin{lstlisting}
def ExtendedLP.IsFeasible [Fintype J] (P : ExtendedLP I J F) : Prop :=
  ∃ p : F∞, P.Reaches p ∧ p ≠ ⊤
\end{lstlisting}
$P$ is bounded by a value $r$ (from below -- we always minimize) if and only if
$P$ reaches only values greater or equal to $r$:
\begin{lstlisting}
def ExtendedLP.IsBoundedBy [Fintype J] (P : ExtendedLP I J F) (r : F) : Prop :=
  ∀ p : F∞, P.Reaches p → r ≤ p
\end{lstlisting}
$P$ is unbounded if and only if $P$ has no finite lower bound:
\begin{lstlisting}
def ExtendedLP.IsUnbounded [Fintype J] (P : ExtendedLP I J F) : Prop :=
  ¬∃ r : F, P.IsBoundedBy r
\end{lstlisting}
Valid extended linear programs are defined as follows (the six conditions were
introduced in \sekt{extensions}):
\newpage
\begin{lstlisting}
structure ValidELP (I J F : Type*) [LinearOrderedField F] extends ExtendedLP I J F where
  hAi : ¬∃ i : I, (∃ j : J, A i j = ⊥) ∧ (∃ j : J, A i j = ⊤)
  hAj : ¬∃ j : J, (∃ i : I, A i j = ⊥) ∧ (∃ i : I, A i j = ⊤)
  hbA : ¬∃ i : I, (∃ j : J, A i j = ⊥) ∧ b i = ⊥
  hcA : ¬∃ j : J, (∃ i : I, A i j = ⊤) ∧ c j = ⊥
  hAb : ¬∃ i : I, (∃ j : J, A i j = ⊤) ∧ b i = ⊤
  hAc : ¬∃ j : J, (∃ i : I, A i j = ⊥) ∧ c j = ⊤
\end{lstlisting}
The following definition says how linear programs are dualized (first, the abbreviation
\texttt{ExtendedLP.dualize} says that $(A, b, c)$ is mapped to $(-A^T, c, b)$;
then \texttt{ValidELP.dualize} says that the dualization of \texttt{ValidELP}
is inherited from \texttt{ExtendedLP.dualize}; the remaining six lines generate
proofs that our six conditions stay satisfied after dualization,
where \texttt{aeply} is a custom tactic based on \texttt{aesop} \cite{Aesop}
and \texttt{apply}):
\begin{lstlisting}
abbrev ExtendedLP.dualize (P : ExtendedLP I J F) : ExtendedLP J I F :=
  ⟨-P.Aᵀ, P.c, P.b⟩

def ValidELP.dualize (P : ValidELP I J F) : ValidELP J I F where
  toExtendedLP := P.toExtendedLP.dualize
  hAi := by aeply P.hAj
  hAj := by aeply P.hAi
  hbA := by aeply P.hcA
  hcA := by aeply P.hbA
  hAb := by aeply P.hAc
  hAc := by aeply P.hAb
\end{lstlisting}
The definition of optimum is, sadly, very complicated (will be explained below):
\begin{lstlisting}
noncomputable def ExtendedLP.optimum [Fintype J] (P : ExtendedLP I J F) : Option F∞ :=
  if ¬P.IsFeasible then
    some ⊤
  else
    if P.IsUnbounded then
      some ⊥
    else
      if hr : ∃ r : F, P.Reaches (toE r) ∧ P.IsBoundedBy r then
        some (toE hr.choose)
      else
        none
\end{lstlisting}
The type \texttt{Option F$\infty$}, which is implemented as
\texttt{Option (Option (Option F))} after unfolding definitions,
allows the following values:
\begin{itemize}
\item \texttt{none}
\item \texttt{some $\bot$} implemented as \texttt{some none}
\item \texttt{some $\top$} implemented as \texttt{some (some none)}
\item \texttt{some (toE r)} implemented as \texttt{some (some (some r))} for any \texttt{r :~F}
\end{itemize}
We assign the following semantics to \texttt{Option F$\infty$} values:
\begin{itemize}
\item \texttt{none} $\ \dots$ invalid finite value (infimum is not attained)
\item \texttt{some $\bot$} $\ \dots$ feasible unbounded
\item \texttt{some $\top$} $\ \dots$ infeasible
\item \texttt{some (toE r)} $\ \dots$ the minimum is a finite value $r$
\end{itemize}
The definition \texttt{ExtendedLP.optimum} first asks if $P$ is feasible;
if not, it returns \texttt{some $\top$} (i.e., the worst value).
When $P$ is feasible, it asks whether $P$ is unbounded; if yes, it returns \texttt{some $\bot$}
(i.e., the best value).
When $P$ is feasible and bounded,
it asks if there is a finite value $r$ such that $P$ reaches $r$ and, at the same time,
$P$ is bounded by $r$; is so, it returns \texttt{some (toE r)}.
Otherwise, it returns \texttt{none}.
Note that we use the verbs ``ask'' and ``return'' metaphorically;
\texttt{ExtendedLP.optimum} is not a computable function; it is just a mathematical
definition (see the keyword \texttt{noncomputable def} above) you can prove things about.
\pagebreak[2]

Finally, we define what opposite values of the type \texttt{Option F$\infty$} are
(note that for any values except \texttt{none} we directly follow the definition
of negation a.k.a.\ unary minus on the extended linearly ordered fields):
\begin{lstlisting}
def OppositesOpt : Option F∞ → Option F∞ → Prop
| (p : F∞), (q : F∞) => p = -q
| _       , _        => False
\end{lstlisting}
For example, \texttt{OppositesOpt (-3) 3} and \texttt{OppositesOpt 5 (-5)} hold.
The last line says that \texttt{OppositesOpt none none} is false.
We~later~show that \texttt{none} is never the optimum of a valid\footnote{
In principle, we could say that \texttt{none} is never the optimum of any
linear program, i.e., $P^\ast$ is defined for every (extended) linear program $P$,
but we did not bother proving it.} linear program, but we do not know it yet.

\newpage 

\section{Formal statements of selected theorems}
\label{statements}

\subsection{Main corollaries}

\begin{lstlisting}
theorem equalityFarkas {I J F : Type*} [Fintype I] [Fintype J]
    [LinearOrderedField F] (A : Matrix I J F) (b : I → F) :
    (∃ x : J → F, 0 ≤ x ∧ A *ᵥ x = b) ≠ (∃ y : I → F, 0 ≤ Aᵀ *ᵥ y ∧ b ⬝ᵥ y < 0)
\end{lstlisting}
\begin{lstlisting}
theorem inequalityFarkas {I J F : Type*} [DecidableEq I] [Fintype I] [Fintype J]
    [LinearOrderedField F] (A : Matrix I J F) (b : I → F) :
    (∃ x : J → F, 0 ≤ x ∧ A *ᵥ x ≤ b) ≠ (∃ y : I → F, 0 ≤ y ∧ 0 ≤ Aᵀ *ᵥ y ∧ b ⬝ᵥ y < 0)
\end{lstlisting}
\begin{lstlisting}
theorem StandardLP.strongDuality {I J R : Type*} [DecidableEq I] [DecidableEq J]
    [Fintype I] [Fintype J] [LinearOrderedField R]
    (P : StandardLP I J R) (hP : P.IsFeasible ∨ P.dualize.IsFeasible) :
    OppositesOpt P.optimum P.dualize.optimum
\end{lstlisting}

\subsection{Main results}

\begin{lstlisting}
theorem fintypeFarkasBartl {J R V W : Type*} [Fintype J] [LinearOrderedDivisionRing R]
    [LinearOrderedAddCommGroup V] [Module R V] [PosSMulMono R V]
    [AddCommGroup W] [Module R W]
    (A : W →ₗ[R] J → R) (b : W →ₗ[R] V) :
    (∃ x : J → V, 0 ≤ x ∧ ∀ w : W, ∑ j : J, A w j • x j = b w) ≠ (∃ y : W, 0 ≤ A y ∧ b y < 0)
\end{lstlisting}

Note that \texttt{→$_l$[R]} encodes a linear map as a special case of a semi-linear map \cite{Semilinear}.

\begin{lstlisting}
theorem extendedFarkas {I J F : Type*} [DecidableEq I]
    [Fintype I] [Fintype J] [LinearOrderedField F]
    (A : Matrix I J F∞) (b : I → F∞)
    (hAi : ¬∃ i : I, (∃ j : J, A i j = ⊥) ∧ (∃ j : J, A i j = ⊤))
    (hAj : ¬∃ j : J, (∃ i : I, A i j = ⊥) ∧ (∃ i : I, A i j = ⊤))
    (hAb : ¬∃ i : I, (∃ j : J, A i j = ⊤) ∧ b i = ⊤)
    (hbA : ¬∃ i : I, (∃ j : J, A i j = ⊥) ∧ b i = ⊥) :
    (∃ x : J → F≥0, A ₘ* x ≤ b) ≠ (∃ y : I → F≥0, -Aᵀ ₘ* y ≤ 0 ∧ b ᵥ⬝ y < 0)
\end{lstlisting}
\begin{lstlisting}
theorem ValidELP.strongDuality {I J F : Type*} [DecidableEq I] [DecidableEq J]
    [Fintype I] [Fintype J] [LinearOrderedField F]
    (P : ValidELP I J F) (hP : P.IsFeasible ∨ P.dualize.IsFeasible) :
    OppositesOpt P.optimum P.dualize.optimum
\end{lstlisting}

\section {Proving the Farkas-Bartl theorem}
\label{bartl}

We prove \ref{finFarkasBartl} and, in the end, we obtain \ref{fintypeFarkasBartl} as corollary.

\begin{theorem}[\customlabel{finFarkasBartl}{\texttt{finFarkasBartl}}]
Let $n$ be a natural number.
Let $R$ be a linearly ordered division ring.
Let $W$ be an $R$-module.
Let $V$ be a linearly ordered $R$-module whose ordering satisfies
monotonicity of scalar multiplication by nonnegative elements on the left.
Let $A$ be an $R$-linear map from $W$ to $(\fin{n} \r R)$.
Let $b$ be an $R$-linear map from $W$ to $V$.
Exactly one of the following exists:
\begin{itemize}
\item nonnegative vector family $x : \fin{n} \r V$ such that, for all $w : W$, we have
$ \sum_{j : \fin{n}}\; (A~w)_j \cdot x_j = b~w $
\item vector $y : W$ such that $A~y \ge 0$ and $b~y < 0$
\end{itemize}
\end{theorem}
The only difference from \ref{fintypeFarkasBartl} is that
\ref{finFarkasBartl} uses $\fin{n} = \{ 0, \dots, n\!-\!1 \}$
instead of an arbitrary (unordered) finite type $J$.

\medskip \noindent
Proof idea:
We first prove that both cannot exist at the same time.
Assume we have $x$ and $y$ of said properties.
We plug $y$ for $w$ and obtain
$ \sum_{j : \fin{n}}\; (A~y)_j \cdot x_j = b~y $.
On the left-hand side, we have a sum of nonnegative vectors,
which contradicts $b~y < 0$.
\smallskip

We prove ``at least one exists'' by induction on $n$.
If $n=0$ then $A~y \ge 0$ is a tautology.
We consider $b$. Either $b$ maps everything to the
zero vector, which allows $x$ to be the empty vector family,
or some $w$ gets mapped to a nonzero vector, where
we choose $y$ to be either $w$ or $(-w)$.
Since $V$ is linearly ordered, one of them satisfies $b~y<0$.
Now we precisely state how the induction step will be.

\begin{lemma}[\customlabel{industepFarkasBartl}{\texttt{industepFarkasBartl}}]
Let $m$ be a natural number.
Let $R$ be a linearly ordered division ring.
Let $W$ be an $R$-module.
Let $V$ be a linearly ordered $R$-module whose ordering satisfies
monotonicity of scalar multiplication by nonnegative elements on the left.
Assume (induction hypothesis) that
for all $R$-linear maps $\bar{A} : W \r (\fin{m} \r R)$
and $\bar{b} : W \r V$, the formula
``$\forall \bar{y} : W \st \bar{A}~\bar{y} \ge 0 \implies \bar{b}~\bar{y} \ge 0$''
implies existence of a nonnegative vector family $\bar{x} : \fin{m} \r V$ such that,
for all $\bar{w} : W$, $ \sum_{i : \fin{m}}\; (\bar{A}~\bar{w})_i \cdot \bar{x}_i = \bar{b}~\bar{w} $.
Let $A$ be an $R$-linear map from $W$ to $(\fin{m\!+\!1} \r R)$.
Let $b$ be an $R$-linear map from $W$ to~$V$.
Assume that, for all $y : W$, $\;A~y \ge 0$ implies $b~y \ge 0$.
We claim there exists a nonnegative vector family $x : \fin{m\!+\!1} \r V$
such that, for all $w : W$, we have
$ \sum_{i : \fin{m+1}}\; (A~w)_i \cdot x_i = b~w $.
\end{lemma}

\noindent
Proof idea:
Let $A_{<m}$ denote a function that maps $(w : W)$ to
$(A~w) \big|_{\fin{m}}$, i.e., $A_{<m}$ is an $R$-linear map
from $W$ to $(\fin{m} \r R)$ that behaves exactly like $A$ where
it is defined.
We distinguish two cases. If, for all $y : W$, the inequality
$A_{<m}~y \ge 0$ implies $b~y \ge 0$, then plug $A_{<m}$
for $\bar{A}$, obtain $\bar{x}$, and construct a vector family $x$ such that
$x_m = 0$ and otherwise $x$ copies $\bar{x}$. We easily check that
$x$ is nonnegative and that
$ \sum_{i : \fin{m+1}}\; (A~w)_i \cdot x_i = b~w $ holds. \pagebreak[2]

In the second case, we have $y'$ such that $A_{<m}~y' \ge 0$
holds but $b~y' < 0$ also holds. We realize that $(A~y')_m < 0$.
We now declare $y := ((A~y')_m)^{-1} \cdot y'$ and observe
$(A~y)_m = 1$. We establish the following facts (proofs are omitted):
\begin{itemize}
\item for all $w : W$, we have $A~(w - ((A~w)_m \cdot y)) = 0$
\item for all $w : W$, the inequality $A_{<m}~(w - ((A~w)_m \cdot y)) \ge 0$
implies $b~(w - ((A~w)_m \cdot y)) \ge 0$
\item for all $w : W$, the inequality $(A_{<m} - (z \mapsto (A~z)_m \cdot (A_{<m}~y)))~w \ge 0$
implies\\ $(b - (z \mapsto (A~z)_m \cdot (b~y)))~w \ge 0$
\end{itemize}
We observe that
$\bar{A} := A_{<m} - (z \mapsto (A~z)_m \cdot (A_{<m}~y))$
and
$\bar{b} := b - (z \mapsto (A~z)_m \cdot (b~y))$
are $R$-linear maps.
Thanks to the last fact, we can apply induction hypothesis to $\bar{A}$ and $\bar{b}$.
We obtain a nonnegative vector family $x'$ such that,
for all $\bar{w} : W$, $ \sum_{i : \fin{m}} (\bar{A}~\bar{w})_i \cdot x'_i = \bar{b}~\bar{w} $.
It remains to construct a nonnegative vector family $x : \fin{m\!+\!1} \r V$
such that, for all $w : W$, we have
$ \sum_{i : \fin{m+1}}\; (A~w)_i \cdot x_i = b~w $.
We choose $x_m = b~y - \sum_{i : \fin{m}} (A_{<m}~y)_i \cdot x'_i$
and otherwise $x$ copies $x'$. We check that our $x$ has the required
properties. \qed

\medskip
We complete the proof of \ref{finFarkasBartl} by applying \ref{industepFarkasBartl}
to $A_{\le n}$ and $b$. Finally, we obtain \ref{fintypeFarkasBartl} from
\ref{finFarkasBartl} using some boring mechanisms regarding equivalence between
finite types.

\section {Proving the Extended Farkas theorem}
\label{extended}

We prove \ref{inequalityFarkas} by applying \ref{equalityFarkas}
to the matrix $(1~|~A)$ where $1$ is the identity matrix of type
$(I \times I) \r F$.

\begin{theorem}[\customlabel{inequalityFarkas-neg}{\texttt{inequalityFarkas\_neg}}]
Let $I$ and $J$ be finite types.
Let $F$ be a linearly ordered field.
Let $A$ be a matrix of type $(I \times J) \r F$.
Let $b$ be a vector of type $I \r F$.
Exactly one of the following exists:
\begin{itemize}
	\item nonnegative vector $x : J \r F$ such that $A \* x \le b$
	\item nonnegative vector $y : I \r F$ such that $(-A^T) \* y \le 0$ and $b \* y < 0$
\end{itemize}
\end{theorem}
Obviously, \ref{inequalityFarkas-neg} is an immediate corollary of \ref{inequalityFarkas}.

\begin{theorem}[\customlabelrestated{extendedFarkas:restated}{\texttt{extendedFarkas}}]
Let $I$ and $J$ be finite types.
Let $F$ be a linearly ordered field.
Let $A$ be a matrix of type $(I \times J) \r F_\infty$.
Let $b$ be a vector of type $I \r F_\infty$.
Assume that $A$ does not have $\bot$ and $\top$ in the same row.
Assume that $A$ does not have $\bot$ and $\top$ in the same column.
Assume that $A$ does not have $\top$ in any row where $b$ has $\top$.
Assume that $A$ does not have $\bot$ in any row where $b$ has~$\bot$.
Exactly one of the following exists:
\begin{itemize}
\item nonnegative vector $x : J \r F$ such that $A \* x \le b$
\item nonnegative vector $y : I \r F$ such that $(-A^T) \* y \le 0$ and $b \* y < 0$
\end{itemize}
\end{theorem}

\noindent
Proof idea: We need to do the following steps in the given order.
\begin{enumerate}
\item Delete all rows of both $A$ and $b$ where $A$ has $\bot$ or $b$ has $\top$
(they are tautologies).
\item Delete all columns of $A$ that contain $\top$
(they force respective variables to be zero).
\item If $b$ contains $\bot$, then $A \* x \le b$ cannot be satisfied,
but $y = 0$ satisfies $(-A^T) \* y \le 0$ and $b \* y < 0$. Stop here.
\item Assume there is no $\bot$ in $b$. Use \ref{inequalityFarkas-neg}.
In either case, extend $x$ or $y$ with zeros on all deleted positions.
\end{enumerate}

\section {Proving the Extended strong LP duality}
\label{dualityELP}

We start with the weak LP duality and then move to the strong LP duality.
We will use \ref{extendedFarkas} in several places.

\begin{lemma}[\customlabel{weakDuality-of-no-bot}{\texttt{weakDuality\_of\_no\_bot}}]
Let $F$ be a linearly ordered field.
Let $P = (A,b,c)$ be a valid linear program over $F_\infty$ such that
neither $b$ nor $c$ contains $\bot$.
If $P$ reaches $p$ and the dual of $P$ reaches $q$,
then $p + q \ge 0$.
\end{lemma}

\noindent
Proof idea:
There is a vector $x$ such that $A \* x \le b$ and $c \* x = p$.
Apply \ref{extendedFarkas} to the following matrix and vector:
$$
\begin{pmatrix}
A \\
c
\end{pmatrix}
\qquad \qquad
\begin{pmatrix}
b \\
p
\end{pmatrix}
$$

\begin{lemma}[\customlabel{no-bot-of-reaches}{\texttt{no\_bot\_of\_reaches}}]
Let $F$ be a linearly ordered field.
Let $P = (A,b,c)$ be a valid linear program over $F_\infty$.
If $P$ reaches any value $p$, then $b$ does not contain $\bot$.
\end{lemma}

\noindent
Proof idea: It would contradict the assumption that
$A$ does not have $\bot$ in any row where $b$ has $\bot$.

\begin{theorem}[\customlabel{weakDuality}{\texttt{weakDuality}}]
Let $F$ be a linearly ordered field.
Let $P$ be a valid linear program over $F_\infty$.
If $P$ reaches $p$ and the dual of $P$ reaches $q$,
then $p + q \ge 0$.
\end{theorem}

\noindent
Proof idea: Apply \ref{no-bot-of-reaches} to P.
Apply \ref{no-bot-of-reaches} to the dual of P.
Finish the proof using \ref{weakDuality-of-no-bot}.

\begin{lemma}[\customlabel{unbounded-of-reaches-le}{\texttt{unbounded\_of\_reaches\_le}}]
Let $F$ be a linearly ordered field.
Let $P$ be a valid linear program over $F_\infty$.
Assume that for each $r$ in $F$ there exists
$p$ in $F_\infty$ such that $P$ reaches $p$
and $p \le r$.
Then $P$ is unbounded.
\end{lemma}

\noindent
Proof idea:
It suffices to prove that for each $r'$ in $F$ there exists
$p'$ in $F_\infty$ such that $P$ reaches $p'$ and $p' < r'$.
Apply the assumption to $r'\!-\!1$.

\begin{lemma}[\customlabel{unbounded-of-feasible-of-neg}{\texttt{unbounded\_of\_feasible\_of\_neg}}]
Let $F$ be a linearly ordered field.
Let $P$ be a valid linear program over $F_\infty$ that is feasible.
Let $x_0$ be a nonnegative vector such that $c \cdot x_0 < 0$
and $A \cdot x_0 + 0 \cdot (-b) \le 0$.
Then $P$ is unbounded.
\end{lemma}

\noindent
Proof idea:
There is a nonnegative vector $x_p$ such that $A \* x_p \le b$ and
$c \* x_p = e$ for some $e \neq \top$.
We apply \ref{unbounded-of-reaches-le}.
In case $e \le r$, we use $x_p$ and we are done.
Otherwise, consider what $c \cdot x_0$ equals to.
If $c \cdot x_0 = \bot$ then we are done.
We cannot have $c \cdot x_0 = \top$ because $c \cdot x_0 < 0$.
Hence $c \cdot x_0 = d$ for some $d$ in $F$.
Observe that the fraction $\frac{r-e}{d}$ is
well defined and it is positive.
Use $x_p + \frac{r-e}{d} \cdot x_0$.

\begin{lemma}[\customlabel{unbounded-of-feasible-of-infeasible}{\texttt{unbounded\_of\_feasible\_of\_infeasible}}]
Let $P$ be a valid linear program such that
$P$ is feasible but the dual of $P$ is not feasible.
Then $P$ is unbounded.
\end{lemma}

\noindent
Proof idea:
Apply \ref{extendedFarkas} to the matrix
$( A \big|_{\{ i : I \,|\, b_i \neq \top \}} )^T$
and the vector $c$.

\begin{lemma}[\customlabel{infeasible-of-unbounded}{\texttt{infeasible\_of\_unbounded}}]
If a valid linear program $P$ is unbounded,
the dual of $P$ cannot be feasible.
\end{lemma}

\noindent
Proof idea:
Assume that $P$ is unbounded, but the dual of $P$ is feasible.
Obtain contradiction using \ref{weakDuality}.

\begin{lemma}[\customlabel{dualize-dualize}{\texttt{dualize\_dualize}}]
Let $P$ be a valid linear program. The dual of the dual of $P$ is exactly $P$.
\end{lemma}

\noindent
Proof idea:
$ -(-A^T)^T = A $

\begin{lemma}[\customlabel{strongDuality-aux}{\texttt{strongDuality\_aux}}]
Let $P$ be a valid linear program such that
$P$ is feasible and the dual of $P$ is also feasible.
There is a value $p$ reached by $P$ and
a value $q$ reached by the dual of $P$ such
that $p + q \le 0$.
\end{lemma}

\noindent
Proof idea:
Apply \ref{extendedFarkas} to the following matrix and vector:
$$
\begin{pmatrix}
	A & 0 \\
	0 & -A^T \\
	c & b
\end{pmatrix}
\qquad \qquad
\begin{pmatrix}
	b \\
	c \\
	0
\end{pmatrix}
$$
In the first case, we obtain a nonnegative vectors $x$ and $y$
such that $A \* x \le b$, $-A^T \* y \le c$, and
$c \* x + b \* y \le 0$.
We immediately see that $P$ reaches $c \* x$,
that the dual of $P$ reaches $b \* y$,
and that the desired inequality holds.

In the second case, we obtain nonnegative vectors $y$ and $x$ and
a nonnegative scalar $z$ such that
$-\!A^T \* y + z \* (-c) \le 0$, $A \* x + z \* (-b) \le 0$,
and $b \* y + c \* x < 0$.
If $z > 0$, we finish the proof using $p := z^{-1} \* c \* x$
and $q := z^{-1} \* b \* y$ (and we do not care that this case
cannot happen in reality).
It remains to prove that $z=0$ cannot happen.
At least one of $b \* y$ or $c \* x$ must be strictly negative.
If $c \* x < 0$, we apply \ref{unbounded-of-feasible-of-neg} to
$P$, which (using \ref{infeasible-of-unbounded}) contradicts
that the dual of $P$ is feasible.
If $b \* y < 0$, we apply \ref{unbounded-of-feasible-of-neg} to
the dual of $P$, which (using \ref{infeasible-of-unbounded} and
\ref{dualize-dualize}) contradicts that $P$ is feasible.

\begin{lemma}[\customlabel{strongDuality-of-both-feasible}{\texttt{strongDuality\_of\_both\_feasible}}]
Let $P$ be a valid linear program such that
$P$ is feasible and the dual of $P$ is also feasible.
There is a finite value $r$ such that
$P$ reaches $-r$ and the dual of $P$ reaches $r$.
\end{lemma}

\noindent
Proof idea:
From \ref{strongDuality-aux} we have a value $p$ reached by $P$ and
a value $q$ reached by the dual of $P$ such that $p + q \le 0$.
We apply \ref{weakDuality} to $p$ and $q$ to obtain $p + q \ge 0$.
We set $r := q$.

\begin{lemma}[\customlabel{optimum-unique}{\texttt{optimum\_unique}}]
Let $P$ be a valid linear program.
Let $r$ be a value reached by $P$ such that $P$ is bounded by $r$.
Let $s$ be a value reached by $P$ such that $P$ is bounded by $s$.
Then $r = s$.
\end{lemma}

\noindent
Proof idea:
We prove $r \le s$ by applying ``$P$ is bounded by $r$'' to
``$s$ is reached by $P$''.
We prove $s \le r$ by applying ``$P$ is bounded by~$s$'' to
``$r$ is reached by $P$''.

\begin{lemma}[\customlabel{optimum-eq-of-reaches-bounded}{\texttt{optimum\_eq\_of\_reaches\_bounded}}]
Let $P$ be a valid linear program.
Let $r$ be a value reached by $P$ such that $P$ is bounded by $r$.
Then the optimum of $P$ is $r$.
\end{lemma}

\noindent
Proof idea:
Apply the axiom of choice to the definition of optimum and use \ref{optimum-unique}.

\begin{lemma}[\customlabel{strongDuality-of-prim-feas}{\texttt{strongDuality\_of\_prim\_feasible}}]
Let $P$ be a valid linear program that is feasible.
Then the optimum of $P$ and the optimum of dual of $P$ are opposites.
\end{lemma}

\noindent
Proof idea:
If the dual of $P$ is feasible as well, use \ref{strongDuality-of-both-feasible}
to obtain $r$ such that $P$ reaches $-r$ and the dual of $P$ reaches $r$.
Using \ref{weakDuality} together with \ref{optimum-eq-of-reaches-bounded},
conclude that the optimum of $P$ is $-r$.
Using \ref{weakDuality} together with \ref{optimum-eq-of-reaches-bounded},
conclude that the optimum of the dual of $P$ is $r$.
Observe that $-r$ and $r$ are opposites.

If the dual of $P$ is infeasible, the optimum of the dual of $P$ is
$\top$ by definition. Using \ref{unbounded-of-feasible-of-infeasible}
we get that the optimum of $P$ is $\bot$.
Observe that $\top$ and $\bot$ are opposites.

\begin{theorem}[\customlabel{optimum-neq-none}{\texttt{optimum\_neq\_none}}]
Every valid linear program has optimum.
\end{theorem}

\noindent
Proof idea:
If a valid linear program $P$ is feasible, the existence of optimum
follows from \ref{strongDuality-of-prim-feas}.
Otherwise, the optimum of $P$ is $\top$ by definition.

\begin{lemma}[\customlabel{strongDuality-of-dual-feas}{\texttt{strongDuality\_of\_dual\_feasible}}]
Let $P$ be a valid linear program whose dual is feasible.
Then the optimum of $P$ and the optimum of dual of $P$
are opposites.
\end{lemma}

\noindent
Proof idea:
Apply \ref{strongDuality-of-prim-feas} to the dual of $P$ and
finish the proof using \ref{dualize-dualize}.

\begin{theorem}[\customlabelrestated{strongDuality:restated}{\texttt{strongDuality}}]
Let $F$ be a linearly ordered field.
Let $P$ be a valid linear program over $F_\infty$.
If $P$ or its dual is feasible (at least one of them),
then the optimum of $P$ and the optimum of dual of $P$
are opposites.
\end{theorem}

\noindent
Proof idea:
Use \ref{strongDuality-of-prim-feas} or \ref{strongDuality-of-dual-feas}.

\section{Counterexamples}
\label{counterexamples}

\subsection{Counterexamples for~\ref{extendedFarkas}}
\label{sec:extendedFarkas:counterexamples}

Recall that~\ref{extendedFarkas} has four preconditions on matrix $A$ and vector $b$.
The following examples show that omitting any of these preconditions makes the theorem
false.
\begin{itemize}
\item If $A$ has $\bot$ and $\top$ in the same row, it may happen that both $x$ and $y$ exist:
$$
A =
\begin{pmatrix}
	\bot & \top \\
	0 & -1
\end{pmatrix}
\qquad \qquad
b = \begin{pmatrix} 0 \\ -1 \end{pmatrix}
\qquad \qquad
x = \begin{pmatrix} 1 \\ 1 \end{pmatrix}
\qquad \qquad
y = \begin{pmatrix} 0 \\ 1 \end{pmatrix}
$$
\item If $A$ has $\bot$ and $\top$ in the same column, it may happen that both $x$ and $y$ exist:
$$
A = \begin{pmatrix} \bot \\ \top \end{pmatrix}
\qquad \qquad
b = \begin{pmatrix} -1 \\ 0 \end{pmatrix}
\qquad \qquad
x = \begin{pmatrix} 0 \end{pmatrix}
\qquad \qquad
y = \begin{pmatrix} 1 \\ 1 \end{pmatrix}
$$
\item If $A$ has $\top$ in a row where $b$ has $\top$, it may happen that both $x$ and $y$ exist:
$$
A = \begin{pmatrix} \top \\ -1 \end{pmatrix}
\qquad \qquad
b = \begin{pmatrix} \top \\ -1 \end{pmatrix}
\qquad \qquad
x = \begin{pmatrix} 1 \end{pmatrix}
\qquad \qquad
y = \begin{pmatrix} 0 \\ 1 \end{pmatrix}
$$
\item If $A$ has $\bot$ in a row where $b$ has $\bot$, it may happen that both $x$ and $y$ exist:
$$
A = \begin{pmatrix} \bot \end{pmatrix}
\qquad \qquad
b = \begin{pmatrix} \bot \end{pmatrix}
\qquad \qquad
x = \begin{pmatrix} 1 \end{pmatrix}
\qquad \qquad
y = \begin{pmatrix} 0 \end{pmatrix}
$$
\end{itemize}

We also claimed in \sekt{extensions} that changing condition $(-A^T) \* y \le 0$
to  $A^T \* y \ge 0$ in  \ref{extendedFarkas} would cause the theorem
to fail even when $A$ has only a single $\bot$ entry. The counterexample is as follows:
$$
\begin{array}{c@{\hspace{30pt}}c}
A =
\begin{pmatrix}
	\bot  \\
	0 
\end{pmatrix}
&
b = \begin{pmatrix} 0 \\ -1 \end{pmatrix}
\end{array}
$$
System $A \* x\le b$ does not have a solution, and neither does system
$A^T \* y \ge 0$, $b \* y < 0$ (over nonnegative vectors $x,y$).

\subsection{Counterexamples for~\ref{ExtendedLP.strongDuality}}\label{sec:ExtendedLP.strongDuality:counterexamples}
Recall that \ref{ExtendedLP.strongDuality} is formulated for {\em valid} LPs,
and the definition of a valid LP has six conditions.
In this section we show that omitting any of these six conditions
makes the theorem false.

\noindent 
Let us consider the following six LPs over $F_\infty$;
all of them are written in the format $(A,b,c)$.
$$
\begin{array}{c@{\hspace{100pt}}c}
P_1 = \left(
	\begin{pmatrix} \bot \\ \top \end{pmatrix},
	\begin{pmatrix} -1 \\ 0 \end{pmatrix},
	\begin{pmatrix} 0 \end{pmatrix}
\right)
&
D_1 = \left(
	\begin{pmatrix} \top & \bot \end{pmatrix},
	\begin{pmatrix} 0 \end{pmatrix},
	\begin{pmatrix} -1 \\ 0 \end{pmatrix}
\right)
\vspace{5pt} \\
P_2=\left(
	\begin{pmatrix} \bot \end{pmatrix}, 
	\begin{pmatrix} \bot \end{pmatrix}, 
	\begin{pmatrix} 0 \end{pmatrix}
\right)
&
D_2=\left(
	\begin{pmatrix} \top \end{pmatrix},
	\begin{pmatrix} 0 \end{pmatrix},
	\begin{pmatrix} \bot \end{pmatrix}
\right)
\vspace{5pt} \\
P_3=\left(
	\begin{pmatrix} \top \\ -1 \end{pmatrix},
	\begin{pmatrix} \top \\ -1 \end{pmatrix},
	\begin{pmatrix} 0 \end{pmatrix}
\right)
&
D_3=\left(
	\begin{pmatrix} \bot & 1 \end{pmatrix},
	\begin{pmatrix} 0 \end{pmatrix},
	\begin{pmatrix} \top \\ -1 \end{pmatrix}
\right)
\end{array}
$$
It can be checked that $D_i$ is the dual of $P_i$ for $i=1,2,3$.
Furthermore, strong duality fails in all cases, since
the optimum of $P_i$ is $0$ and the optimum of $D_i$ is $\bot$ for $i=1,2,3$.
Each of $P_1,D_1,P_2,D_2,P_3,D_3$ violates exactly one of the conditions
in the definition of a valid LP.

\section{Related work}

\subsection{Farkas-like theorems}

There is a substantial body of work on linear inequalities and
linear programming formalized in Isabelle.
While our work is focused only on proving mathematical theorems,
the work in Isabelle is motivated by the development of SMT solvers.

\begin{itemize}
\item 
Bottesch, Haslbeck, Thiemann \cite{Farkas-AFP} proved a variant of
\ref{equalityFarkas} for $\delta$-rationals
by analyzing a specific implementation of the Simplex algorithm
\cite{Simplex-AFP} by Marić, Spasić, Thiemann.
\item
Bottesch, Raynaud, Thiemann \cite{Linear-AFP} proved
the Fundamental theorem of linear inequalities as well
as both \ref{equalityFarkas} and \ref{inequalityFarkas}
for all linearly ordered fields, alongside with
the Carathéodory’s theorem and the Farkas-Minkowski-Weyl theorem.
They also investigated systems of linear mixed-integer inequalities.
\item 
Thiemann himself \cite{Duality-AFP} then proved the strong LP duality
in the asymmetric version \eqref{eq:StandardLP.strongDuality:assym}.
\end{itemize}

Sakaguchi\footnote{\url{https://github.com/pi8027/vass}} proved a version of \ref{equalityFarkas}
for linearly ordered fields in Rocq, using the Fourier-Motzkin elimination.
Allamigeon and Katz \cite{Simplex-Coq} made a large contribution to
the study of convex polyhedra in Rocq\:---\:among
other results, they proved a version of \ref{equalityFarkas}
for linearly ordered fields as well as the strong LP duality
in the asymmetric version \eqref{eq:StandardLP.strongDuality:assym}.

\subsection{Hahn-Banach theorems}

Several Hahn-Banach theorems have been formalized in Lean as a part of
Mathlib.
In particular, it would have been possible to prove \ref{equalityFarkas}
for reals using the Hahn-Banach separation theorem\footnote{\url{https://github.com/leanprover-community/mathlib4/blob/master/Mathlib/Analysis/LocallyConvex/Separation.lean}} for a convex closed set
and a point.
It would then have been possible to extend the result to other
duality theorems for reals.
To our knowledge, this theorem however cannot be used
to prove \ref{equalityFarkas} for an arbitrary linearly ordered field.
Therefore, we decided to prove \ref{equalityFarkas} without appealing to
the geometry and topology.

Note that certain Hahn-Banach theorems have also
been formalized in Mizar \cite{HahnBanach-Mizar},
Alf \cite{HahnBanach-Alf}, Isabelle \cite{HahnBanach-Isabelle},
and Rocq \cite{HahnBanach-Coq}.

\section {Conclusion}

We formally verified several Farkas-like theorems in Lean 4.
We extended the existing theory to a new setting where some
coefficient can carry infinite values. We realized that the
abstract work with modules over linearly ordered division rings
and linear maps between them was fairly easy to carry on in
Lean 4 thanks to the library Mathlib that is perfectly suited
for such tasks. In contrast, manipulation with matrices got
tiresome whenever we needed a not-fully-standard operation.
It turns out Lean 4 cannot automate case analyses unless they
take place in the ``outer layers'' of formulas. Summation
over subtypes and summation of conditional expression made
us developed a lot of ad-hoc machinery which we would have
preferred to be handled by existing tactics. Another area
where Lean 4 is not yet helpful is the search for counterexamples.
Despite these difficulties, we find Lean 4 to be an excellent
tool for elegant expressions and organization of
mathematical theorems and for proving them formally.

\section*{Acknowledgements}

We would like to thank David Bartl and Jasmin Blanchette for
frequent consultations.
We would also like to express gratitude
to Henrik Böving for a help with generalization from
extended rationals to extended linearly ordered fields
and to Andrew Yang for
the proof of \texttt{Finset.univ\_sum\_of\_zero\_when\_not}.
We would also like to acknowledge Antoine Chambert-Loir,
Apurva Nakade, Yaël Dillies, Richard Copley,
Edward van de Meent, Markus Himmel, Mario Carneiro,
and Kevin Buzzard.

\printbibliography

\newpage

\section*{Appendix: dependencies between theorems}

Theorems are in black. Selected lemmas are in gray.
What we consider to be the main theorems are
denoted by blue background.
What we consider to be the main corollaries are
denoted by yellow background.

\begin{center}
\includegraphics[width=0.75\textwidth]{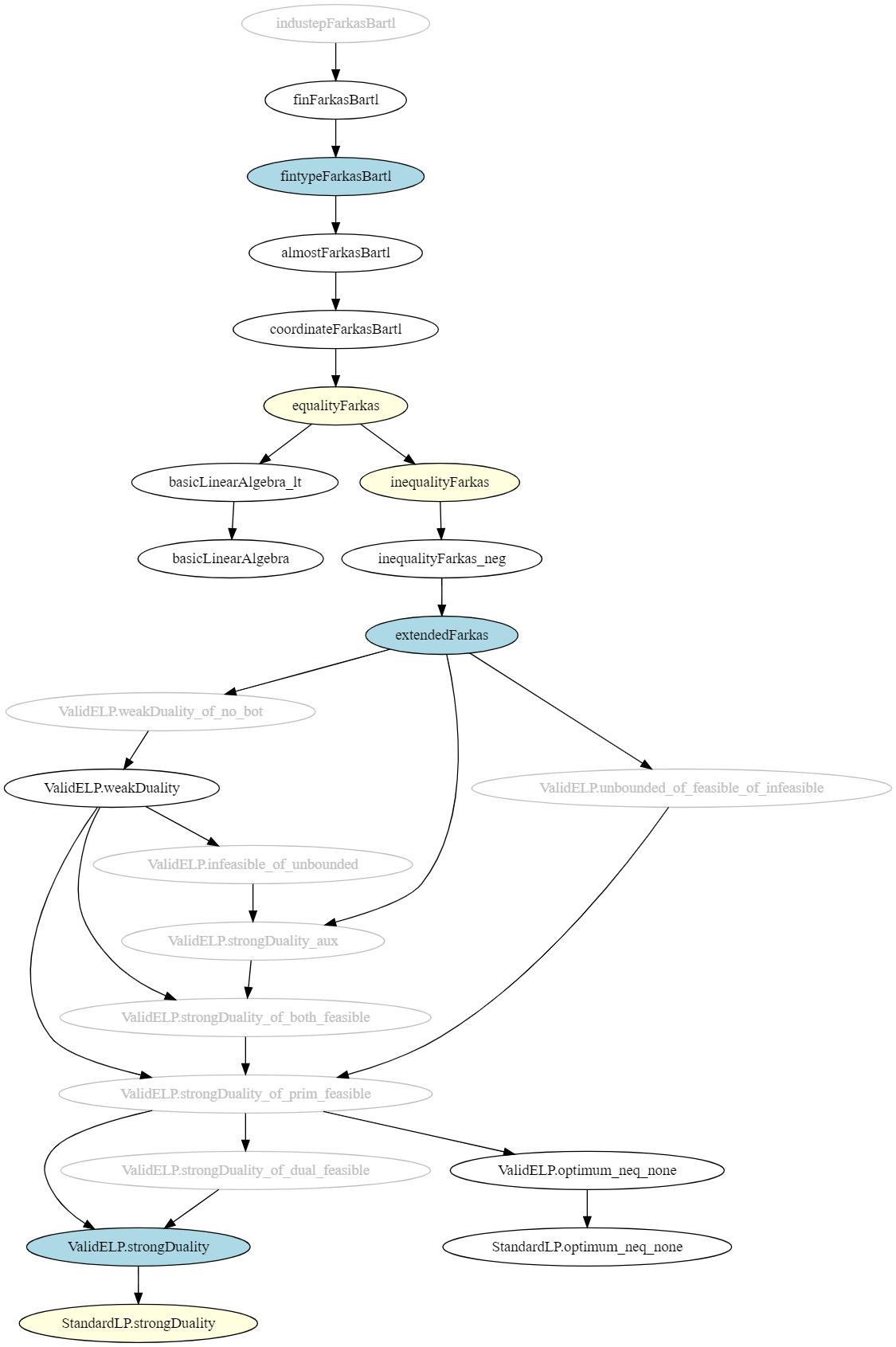}
\end{center}

\end{document}